\documentclass[11pt]{article}

\usepackage{amsfonts}

\begin{document}




\parindent = 0 pt

\newtheorem{thm}{Theorem}[section]

\newtheorem{lem}[thm]{Lemma}

\newtheorem{pre-note}[thm]{Note}
\newenvironment{note}{\begin{pre-note}\rm}{\end{pre-note}}

\newtheorem{pre-proof}[thm]{Proof}
\newenvironment{proof}{\begin{pre-proof}\rm}{$\quad\bullet$\end{pre-proof}}

\newtheorem{pre-example}[thm]{Example}
\newenvironment{example}{\begin{pre-example}\rm}{\end{pre-example}}

\newtheorem{pre-numeq}[thm]{(}
\newenvironment{numeq}{\begin{pre-numeq}\rm\bf)$\quad\displaystyle}{$\end{pre-numeq}}


\newcommand{\RR}{\mathbb{R}}
\newcommand{\CC}{\mathbb{C}}
\newcommand{\HH}{\mathbb{H}}

\newcommand{\Real}{{\rm Re}}

\newcommand{\NB}{{\bf NB}}
\newcommand{\bH}{{\bf H}}
\newcommand{\Un}{{\bf 1}}
\newcommand{\bP}{{\bf P}}

\begin{flushright}
math/9910055\\
msc2000: Primary 11M26, Secondary 30D50\\
\end{flushright}

\begin{center}
{\bf \Large A NOTE ON NYMAN'S EQUIVALENT FORMULATION OF THE RIEMANN HYPOTHESIS}\\
\vskip 1 cm
{\Large Jean-Fran\c{c}ois Burnol}\\
\vskip 0.5 cm
October 1999 (v1), January 2001 (v2)\\
\vskip 0.5cm
\end{center}


{\bf Abstract:} A certain subspace of $L^2((0,1), dt)$ has been considered by
Nyman, Beurling, and others, with the result that the constant function {\bf 1}
belongs to it if and only if the Riemann Hypothesis holds. I show in this note
that the product $\prod_{\zeta(\rho)=0,
\Real(\rho)>{1\over2}}\left|{{1-\rho}\over\rho}\right|$ is the norm of the
projection of {\bf 1} to this subspace. This provides a quantitative refinement
to Nyman's theorem.






\section{Nyman's criterion for the Riemann Hypothesis}

We first recall a theorem of Nyman \cite{Nyman}. We write $\{u\}$ for the
fractional part of the real number $u$. To each $\alpha$ in $]0,1[$ we associate
the function $\rho_\alpha$ of the positive real number $t$ defined by
$\rho_\alpha(t) = \left\{\alpha\over t\right\} - \alpha\left\{1\over
t\right\}$. Of course $\rho_\alpha$ vanishes identically on $]1,+\infty[$, and
is bounded on $]0,1]$. We let $\NB$ be the closed subspace of $\bH = L^2((0,1),
dt)$ spanned by the functions $\rho_\alpha$ for $0<\alpha<1$. We consider both
$\NB$ and $\bH$ as closed subspaces of $L^2((0,\infty), dt)$. We let $\Un$ be
the step function with value 1 for $0<t\leq 1$, and 0 for $t>1$. We note that
both $\bH$ and $\NB$ are invariant under the semi-group of unitary contractions
$U(\lambda):f(t)\mapsto {1\over\sqrt{\lambda}}f({t\over\lambda})$, $0<\lambda
\leq 1$. Indeed $U(\lambda)\cdot\rho_\alpha =
{1\over\sqrt{\lambda}}(\rho_{\alpha\lambda} - \alpha\rho_\lambda)$. So $\NB$
contains the constant function $\Un$ if and only if it contains all step
functions if and only if it actually coincides with all of $\bH$.

\begin{thm}[Nyman, 1950 \cite{Nyman}]\label{Nyman}
The Riemann
Hypothesis holds if and only if the function $\Un$ belongs to $\NB$ (hence if and only
if $\NB = \bH$.)
\end{thm}

The paper by Balazard and Saias (\cite{Balazard-Saias}) is an exposition of the
use of the theory of Hardy spaces and especially of fundamental results on
invariant subspaces by Beurling (\cite{BeurlingActa}) and Lax (\cite{Lax}) in
this context. Nyman's criterion has stimulated  various developments in recent
years, theoretical as well as numerical (see \cite{Baez}, \cite{BaBaLaSa} and
references therein).  We will work exclusively in an Hilbert space setting but
should not fail to mention the extension by Beurling (\cite{Beurling}) of
Nyman's work to an $L^p$-setting, leading to a  characterization
of the absence of zeros with real part (strictly) greater than $1\over
p$. Bercovici and Foias (\cite{Bercovici-Foias}, \cite{Ber98}) have shown how
the geometry of the contraction semi--group $U(\lambda),\ \lambda \leq 1$  can
also be the basis of a proof in the $L^p$--case.  We prove:

\begin{thm}\label{MonTheoreme}
Let $\bP$ be orthogonal projection onto $\NB$. We have (counting the zeros with
their multiplicities)
$$\Vert \bP({\bf 1})\Vert = \prod_{\zeta(\rho) = 0,\
\Real(\rho)>{1\over2}}\left|{{1-\rho}\over\rho}\right|$$
\end{thm}

Clearly this implies Nyman's criterion. The proof requires some classical
results of harmonic analysis (Hardy spaces, the factorization theorem, the
Beurling--Lax description of invariant subspaces) for which accessible reviews
are given in the books by Dym--McKean (\cite{Dym-McKean}) and Hoffman
(\cite{Hoffman}). Specifically we will put to use an old general argument of
Beurling (\cite[Theorem II]{BeurlingActa}) (originating in the study of the
Hardy space of the unit disc). It is natural to discuss {\bf \ref{MonTheoreme}}
together with a proof (whose idea is given in \cite{Burnol}) of the elegant
following result:

\begin{thm}[Balazard-Saias-Yor, \cite{Balazard-Saias-Yor}]\label{Bal}One has 
$${1\over2\pi}\int_{\Real(w) = {1\over2}} {\log|\zeta(w)|\over|w|^2}\,|dw| =
\sum_{\zeta(\rho) = 0,\ \Real(\rho)>{1\over2}}\log\left|{\rho\over
1-\rho}\right|$$\end{thm}


\section{Proofs of \ref{Nyman}, \ref{MonTheoreme} and \ref{Bal}}

As in the papers cited above on this subject the main tool is the Fourier
Transform, here in its multiplicative version $f(t) \in L^2((0,\infty), dt)\
\mapsto \widehat{f}(s) = \int_0^\infty f(t)\,t^{s-1}\,dt$.  The integral is to
be understood in the $L^2$--sense, with $\Real(s) = {1\over2}$. Indeed we really
want to look at $\sqrt{t}f(t)$ in $L^2((0,\infty), {dt\over t})$ and at its
transform in the dual group $\int_0^\infty \sqrt{t}f(t)\,t^{i\tau}\,{dt\over
t}$, $\tau \in \RR$. With $s = {1\over2} + i\tau$ we end up with the formula
above. If $f(t)$ actually belongs to $L^2((0,1), dt)$, then $\widehat{f}(s) =
\int_0^\infty f(t)\,t^{s-1}\,dt$ makes sense as an analytic function in the
half--plane $\Real(s) > {1\over2}$ and the set of such $\widehat{f}$ is
characterized by the Paley--Wiener Theorem (\cite[Chapter 3]{Dym-McKean}) as the
Hardy space $\HH^2$ of analytic functions $h(s)$ whose $L^2$--norms on vertical
lines ${1\over2\pi}\int_{\Real(s) = \sigma>{1\over2}}{\vert h(s) \vert}^2\,
|ds|$ are bounded independently of $\sigma$. Such an analytic function $h(s)$
has (pointwise almost everywhere) non-tangential limits $h({1\over2}+ i\tau)$
which are also obtained as the Fourier--Mellin Transform of $f\in
L^2((0,1), dt)$. Furthermore, there is a Cauchy formula:
\begin{numeq}\label{Cauchy} 
\qquad\Real(s)>{1\over2}\quad\Rightarrow\quad h(s) =
{1\over2\pi}\int_{-\infty}^{+\infty} {h({1\over2} + i\tau)\over s- {1\over2} -
i\tau}\,d\tau
\end{numeq}

We note in passing that $\widehat{\Un}(s) = {1\over s}$. The zeta-function appears in
this story thanks to $$\int_0^1 \left\{1\over t\right\}t^{s-1}\, dt = {1\over
s-1} - {\zeta(s)\over s}$$


Multiplication by $s-1\over s$ leaves invariant $\HH^2$, so ${(s-1)\zeta(s)\over
s^2}$ belongs to $\HH^2$ (\cite{Burnol}; see also \cite{Ehm}). As reviewed in
\cite{Dym-McKean}, \cite{Hoffman} a general fact following from this is the absolute
convergence for any $s$ in the half-plane $\Real(s) > {1\over2}$ of the Blaschke
product $$B(s) = \prod_{\zeta(\rho) = 0,\ \Real(\rho)>{1\over2}}{s - \rho \over
s - (1 -\bar{\rho})}{1 - \bar{\rho}\over\rho}\left|{\rho\over 1-\rho}\right|$$
It is known that such a product built from the zeros of an element of $\HH^2$ is
an \emph{inner function}, that is an analytic function bounded by $1$ in the
half-plane whose non-tangential limits on the critical line have modulus $1$
(almost everywhere). The expression as an infinite product might cease to make
sense (pointwise) for the boundary values, but in the case at hand it is
absolutely convergent for all $s$ in $\CC$, except at the possible poles $1
-\bar{\rho}$. A general fact is Norbert Wiener's  theorem on the gain of a
causal filter: it implies that the integral involved in {\bf \ref{Bal}} is
absolutely convergent. Furthermore there is a factorization formula
(\cite[chap. 8]{Hoffman}):
$${(s-1)\zeta(s)\over s^2} = c\cdot A^{s-{1\over2}}\cdot
\exp\left({1\over2\pi}\int_{\Real(w) = {1\over2}} \log\left|{\zeta(w)\over
w}\right|\,{s +  w - 2sw\over s - w}\,{|dw|\over |w|^2}\right)\cdot B(s)$$ where
$|c|=1$, $0<A\leq 1$. There could have been a further term associated to a
singular measure on the critical line, but the analytic continuation of
$\zeta(s)$ implies its non-existence. A value of $A$ strictly less than 1 would
imply a much too quick decrease of $\zeta(s)$ when $s$ is real and goes to
$+\infty$, so $A = 1$ (see \cite{Bercovici-Foias}, also \cite{Burnol}). Looking
at real values of $s$ we see that $c$ is real and positive so $c=1$. We thus
have:
$${(s-1)\zeta(s)\over s^2} = \exp\left({1\over2\pi}\int_{\Real(w) = {1\over2}}
\log\left|{\zeta(w)\over w}\right|\,{s + w - 2sw\over s - w}\,{|dw|\over
|w|^2}\right)\cdot B(s)$$ We note that the same applies to  ${1\over s}$  and
gives ${1\over s}  = \exp\left( - {1\over 2\pi}\int_{\Real(w) = {1\over2}}
\log|w|\,{s + w - 2sw\over s - w}\,{|dw|\over |w|^2}\right)$ so that we finally
conclude (for $\Real(s)>{1\over2}$):
\begin{numeq}\label{factor}
{s-1\over s}\zeta(s) = \exp\left({1\over2\pi}\int_{\Real(w) = {1\over2}}
\log\left|{\zeta(w)}\right|\,{s + w - 2sw\over s -
w}\,{|dw|\over |w|^2}\right)\cdot B(s)\end{numeq}
from which the Balazard-Saias-Yor formula {\bf \ref{Bal}} follows by
specializing to $s = 1$. One has furthermore
\begin{numeq}\label{rhohat}
\qquad\qquad\widehat{\rho_\alpha}(s) = {{\alpha-\alpha^s}\over s}\zeta(s)
\end{numeq} 
and the Beurling-Lax (see \cite[chap. 7]{Hoffman}) general description of
invariant subspaces of the Hardy spaces then leads to the
conclusion that the Mellin transform of $\NB$ exactly coincides with
$B(s)\HH^2$. Indeed on one hand the ``common denominator'' of the Blaschke
products associated to the functions $\widehat{\rho_\alpha}(s)$ is $B(s)$ and on
the other hand none of these functions has a singular inner factor, nor  a
special inner factor $A^{s-{1\over2}}$.

\begin{lem} The orthogonal projection of ${1\over s}$ on $B(s)\HH^2$ is
  $B(1){B(s)\over s}$.
\end{lem}



\begin{proof}
The argument is essentially taken from Beurling: see \cite[Theorem
II]{BeurlingActa}. We check that ${1\over s} - B(1) {B(s)\over s}$ is
perpendicular to $B(s)h(s)$, for any $h(s) \in \HH^2$ :
\begin{eqnarray*}
&\ & {1\over2\pi}\int_{\Real(w) = {1\over2}} \overline{\left({1\over w} - B(1)
{B(w)\over w}\right)}\,B(w)h(w)\,|dw|\\ &=& {1\over2\pi}\int_{\Real(w) =
{1\over2}} {\left(B(w) - B(1)\right)h(w)\over \overline{w}}\,|dw|\\ &=&
{1\over2\pi}\int_{\Real(w) = {1\over2}} {\left(B(w) - B(1)\right)h(w)\over
1-w}\,|dw|\\ &=& \left(B(1) - B(1)\right)h(1) = 0\\
\end{eqnarray*}
where $\vert B(w)\vert = 1$ on $\Real(w) = {1\over2}$ and then
({\bf\ref{Cauchy}}) were used. 
\end{proof}

Hence the norm of the projection of $\Un$ on $\NB$ is $B(1)$ as stated in
{\bf\ref{MonTheoreme}} (we note that $ {B(s)\over s}$ is of norm 1 as $B(s)$
is of modulus 1 on the critical line.) We can also interpret Theorem
{\bf\ref{MonTheoreme}} as a statement involving a conditional expectation for a 
certain Gaussian system (work in progress.)

\begin{flushleft}
Jean-Fran\c{c}ois
Burnol\\
Universit\'e de Nice-Sophia-Antipolis\\
Laboratoire J.-A.
Dieudonn\'e\\
Parc Valrose\\
F-06108 Nice C\'edex 02\\
France\\
{\tt burnol@math.unice.fr}\\
\end{flushleft}

\clearpage


\begin{thebibliography}{999}

\bibitem{Baez}
         L.~B\'aez-Duarte,
        \emph{``On Beurling's real variable reformulation of the Riemann
          hypothesis''},
         Adv. in Maths. {\bf 101} (1993), 10-30.

\bibitem{BaBaLaSa}
        L.~B\'aez-Duarte, M.~Balazard, B.~Landreau and E.~Saias,
        \emph{``Notes sur la fonction $ \zeta $ de Riemann,~{\bf 3}''\/}
         Adv. in Maths. {\bf 149} (2000), 130-144.

\bibitem{Balazard-Saias}
         M.~Balazard, E.~Saias,
         \emph{``The Nyman--Beurling equivalent form for the Riemann hypothesis''\/},
         Expo. Math. {\bf 18} (2000), 131-138.

\bibitem{Balazard-Saias-Yor}
        M.~Balazard, E.~Saias, M.~Yor,
        \emph{``Notes sur la fonction $\zeta$ de Riemann, {\bf 2}''\/},
        Advances in Math. {\bf 143}, 284-287, (1999).

\bibitem{Bercovici-Foias}
        H.~Bercovici, C.~Foias,
        \emph{``A real variable restatement of Riemann's Hypothesis''\/},
        Israel J. of Math.\ {\bf 48}, 57--68, (1984).


\bibitem{Ber98} H.~Bercovici, C.~Foias, \emph{``On the Zorn spaces in
Beurling's approach to the Riemann hypothesis''},  Analysis and
topology, 143--149,  World Sci. Publishing, River Edge, NJ, 1998.

\bibitem{Beurling}
        A.~Beurling,
        \emph{``A closure problem related to the Riemann Zeta--function''\/},
        Proc.\ Nat.\ Acad.\ Sci.\ {\bf 41}, 312--314, (1955).

\bibitem{BeurlingActa}
         A.~Beurling,
         \emph{``On two problems concerning linear transformations in Hilbert space''\/},
         Acta Mathematica\ {\bf 81}, 239--255, (1949). 

\bibitem{Burnol}
         J.-F.~Burnol,
         \emph{``An adelic causality problem related to abelian L-functions''\/},
        to be published in the Journal of Number theory.



\bibitem{Dym-McKean}
         H.~Dym, H.~P.~McKean,
         \emph{``Fourier Series and Integrals''\/},
         Academic Press (1972).

\bibitem{Ehm}
        W.~Ehm,
        \emph{``A family of probability densities related to the Riemann
zeta function''}, this volume.



\bibitem{Hoffman}
        K.~Hoffman,
        \emph{``Banach spaces of analytic functions''\/},
        Prentice-Hall, Inc. (1962). (now available as a Dover Publication, 1988).

\bibitem{Lax}
         P.~Lax,
         \emph{``Translation invariant subspaces''\/},
         Acta Mathematica\ {\bf 101}, (1959).


\bibitem{Nyman}
        B.~Nyman,
        \emph{``On some groups and semigroups of translations''\/},
        Thesis, Uppsala (1950).


\end{thebibliography}
\end{document}